\documentclass[journal,twoside]{IEEEtran}
\usepackage{cite}

\ifCLASSINFOpdf
\else
\fi
\usepackage[cmex10]{amsmath}
\usepackage{array}

\usepackage{mdwmath}
\usepackage{mdwtab}
\usepackage[caption=false,font=footnotesize]{subfig}
\usepackage{stfloats}
\usepackage{url}

\usepackage[american,fulldiode]{circuitikz} 
\usepackage{graphicx} 
\usepackage{multirow}
\usepackage{amsfonts}
\usepackage{bbm}
\usepackage{algorithm}
\usepackage[noend]{algpseudocode}
\usepackage{multirow}
\usepackage[section]{placeins}
\usepackage{epstopdf}

\usepackage{slashbox}
\usepackage{amssymb}

\usepackage{lipsum,amsmath,multicol}

\usepackage{multirow}

\usepackage{float}
\usepackage{tablefootnote}
\usepackage[multiple]{footmisc}

\makeatletter
\def\ps@headings{%
\def\@oddhead{\mbox{}\scriptsize\rightmark \hfil \thepage}%
\def\@evenhead{\scriptsize\thepage \hfil \leftmark\mbox{}}%
\def\@oddfoot{}%
\def\@evenfoot{}}
\def\BState{\State\hskip-\ALG@thistlm}
\makeatother

\allowdisplaybreaks

\usepackage{varwidth}
\newcolumntype{M}[1]{>{\begin{varwidth}[t]{#1}}l<{\end{varwidth}}}

\algnewcommand\algorithmicforeach{\textbf{for each}}
\algdef{S}[FOR]{ParForEach}[1]{\algorithmicforeach\ #1\ \algorithmicdo ~\textbf{(in parallel)}}
\algdef{S}[FOR]{ForEach}[1]{\algorithmicforeach\ #1\ \algorithmicdo}
\newcommand{\Break}{\State \textbf{break} }

\begin{document}
%
\title{Computing the Feasible Spaces of \\ Optimal Power Flow Problems}
%
%
%
\author{Daniel K. Molzahn, \IEEEmembership{Member, IEEE}
\thanks{Argonne National Laboratory, Energy Systems Div.: dmolzahn@anl.gov}}
\maketitle

\begin{abstract}
The solution to an optimal power flow (OPF) problem provides a minimum cost operating point for an electric power system. The performance of OPF solution techniques strongly depends on the problem's feasible space. This paper presents an algorithm for provably computing the entire feasible spaces of small OPF problems to within a specified discretization tolerance. Specifically, the feasible space is computed by discretizing certain of the OPF problem's inequality constraints to obtain a set of power flow equations. All solutions to the power flow equations at each discretization point are obtained using the Numerical Polynomial Homotopy Continuation (NPHC) algorithm. To improve computational tractability, ``bound tightening'' and ``grid pruning'' algorithms use convex relaxations to eliminate the consideration of discretization points for which the power flow equations are provably infeasible. The proposed algorithm is used to generate the feasible spaces of two small test cases.
\end{abstract}


\begin{IEEEkeywords}
Optimal power flow, Feasible space, Convex optimization, Global solution
\end{IEEEkeywords}

%
\IEEEpeerreviewmaketitle

\section{Introduction}
\label{l:introduction}
%
%
%
%
\IEEEPARstart{O}{ptimal} power flow (OPF) is one of the key problems in power system optimization.  The OPF problem seeks an optimal operating point in terms of a specified objective function (e.g., minimizing generation cost, matching a desired voltage profile, etc.). Equality constraints are dictated by the network physics (i.e., the power flow equations) and inequality constraints are determined by engineering limits on, e.g., voltage magnitudes, line flows, and generator outputs.

The OPF problem is non-convex due to the non-linear power flow equations, may have local optima~\cite{bukhsh_tps}, and is generally NP-Hard~\cite{lavaei_tps,bienstock2015nphard}, even for networks with tree topologies~\cite{NPhard}. Since first being formulated by Carpentier in 1962~\cite{carpentier}, a broad range of solution approaches have been applied to OPF problems, including successive quadratic programs, Lagrangian relaxation, heuristic optimization, and interior point methods~\cite{opf_litreview1993IandII,ferc4}. Many of these approaches are computationally tractable for large OPF problems. However, despite often finding global solutions~\cite{molzahn_lesieutre_demarco-global_optimality_condition}, these approaches may fail to converge or converge to a local optimum~\cite{bukhsh_tps,ferc5}.

Recently, there has been significant effort focused on convex relaxations of the OPF problem. These include relaxations based on semidefinite programming (SDP)~\cite{lavaei_tps,molzahn_holzer_lesieutre_demarco-large_scale_sdp_opf,pscc2014,madani2014,molzahn_hiskens-sparse_moment_opf,ibm_opf,josz_molzahn-complex_hierarchy,cedric_msdp}, second-order cone programming (SOCP)~\cite{low_tutorial,coffrin2015qc,sun2015,ghaddar2015}, and linear programming (LP)~\cite{bienstock2014,coffrin2016pscc}. In contrast to traditional approaches, convex relaxations provide a lower bound on the optimal objective value, can certify problem infeasibility, and, in many cases, provably yield the global optimum.

The performance of both traditional algorithms and convex relaxations strongly depends on the OPF problem's feasible space characteristics. Accordingly, understanding OPF feasible spaces is crucial for algorithmic research. Characterizing the feasible spaces of OPF problems has been an important research topic~\cite{hiskens2001,bernie_opfconvexity,hill2008,bukhsh_tps,madani2014,lavaei_geometry,chandra2015equilibria,bienstock2015nphard,dj2015,NPhard}.
This paper proposes an algorithm for computing the feasible spaces of small OPF problems. Visualizations resulting from the computed feasible spaces increase researchers' understanding of challenging problems and aid in improving solution algorithms.

The feasible spaces of some OPF problems can be computed analytically. For instance, OPF problems for two-bus systems have analytic solutions~\cite{alves2002,bukhsh_tps,pscc2014}. Exploiting problem symmetries enables explicit expressions for the feasible spaces of other problems~\cite{bernie_opfconvexity}. However, analytic solution is limited to a small set of special cases.

Related work focuses on the \emph{feasibility boundary} of the power flow equations (i.e., the set of parameters for which small parameter changes results in insolvability of the power flow equations). There have been many research efforts in computing the distance to the power flow feasibility boundary for voltage collapse studies, e.g.,~\cite{venikov1975,voltage_collapse,cpf,alvarado1994}. These approaches generally provide small regions (often a single point) that are on the boundary of the feasible space of the power flow equations. A more general continuation-based approach is developed in~\cite{hiskens2001}. Starting from a feasible point, the approach in~\cite{hiskens2001} uses a continuation method to find a point on the power flow feasibility boundary. By freeing a single parameter (e.g., active power injection at one bus), the approach in~\cite{hiskens2001} uses continuation to trace curves that lie on the power flow feasibility boundary. The approach in~\cite{hiskens2001} is computationally tractable for large problems. However, it is difficult to certify that the approach in~\cite{hiskens2001} captures the entire feasible space due to certain non-convexities such as disconnected components. Further, the approach in~\cite{hiskens2001} does not consider all inequality constraints relevant to OPF problems.

The algorithm proposed in this paper is guaranteed to compute the entire OPF feasible space (to within a specified discretization tolerance) for small OPF problems. Specifically, the proposed algorithm discretizes certain inequalities in an OPF problem into equality constraints that take the form of power flow equations. The Numerical Polynomial Homotopy Continuation (NPHC) algorithm~\cite{guo1990,mehta2014a,SW:05,BHSW06} is then used to compute all power flow solutions at each discretization point. The guarantees inherent to the NPHC algorithm ensure the capturing of the entire OPF feasible space. The proposed algorithm is similar to that used in the software Paramotopy~\cite{paramotopy} for visualizing the effects of parameter variation in general polynomial systems.

To improve computational tractability, convex relaxations are employed to eliminate the consideration of infeasible discretization points. Specifically, a hierarchy of ``moment'' relaxations is used to tighten the right hand sides of the OPF problem's inequality constraints. A ``grid pruning'' algorithm is then used to eliminate discretization points that are outside the relaxation's feasible space and therefore provably infeasible.

Many industrially relevant OPF problems have thousands to tens-of-thousands of buses. The proposed feasible space computation algorithm is limited to much smaller problems due to the intractability of NPHC for large problems. Fortunately, there are many small OPF problems with interesting feasible spaces. Further, experience with the moment relaxations of OPF problems suggests that many challenges inherent to large problems are related to non-convexities associated with small regions of the large problems~\cite{molzahn_hiskens-sparse_moment_opf}. By enabling detailed studies of small problems, the proposed algorithm provides the basis for future work in characterizing the physical features that give rise to challenging OPF problems.

The main contributions of this paper are twofold: 1)~Proposal of an OPF-specific algorithm that is guaranteed to compute the complete feasible spaces of small problems. This algorithm is particularly relevant for studies of OPF problems that challenge both traditional solvers and convex relaxation approaches. 2)~The use of convex relaxations to eliminate provably infeasible points, thereby significantly improving computational tractability.

This paper is organized as follows. Section~\ref{l:opf} describes the OPF problem. Section~\ref{l:discretization} presents the proposed discretization approach and the NPHC algorithm used to solve the power flow equations at each discretization point. Section~\ref{l:convex} discusses the application of a hierarchy of convex relaxations to eliminate provably infeasible grid points. Section~\ref{l:examples} applies these techniques to two OPF problems: the five- and nine-bus systems in~\cite{bukhsh_tps}. Section~\ref{l:conclusion} concludes the paper.

\section{Overview of the OPF Problem}
\label{l:opf}

This section presents an OPF formulation in terms of complex voltages, active and reactive power injections, and apparent power line flow limits. Consider an $n$-bus power system, where $\mathcal{N} = \left\lbrace 1,\ldots, n \right\rbrace$ is the set of all buses, $\mathcal{G}$ is the set of generator buses, $\mathcal{S}$ is the index of the bus that fixes the angle reference, and $\mathcal{L}$ is the set of all lines. Let $P_{Di} + \mathbf{j} Q_{Di}$ represent the active and reactive load demand at bus~$i\in\mathcal{N}$, where $\mathbf{j}=\sqrt{-1}$. Let $V_i = V_{di} + \mathbf{j}V_{qi}$ represent the complex voltage phasor at bus~$i \in \mathcal{N}$. Superscripts ``max'' and ``min'' denote specified upper and lower limits. Buses without generators have maximum and minimum generation set to zero. Let $\mathbf{Y} = \mathbf{G}+\mathbf{j}\mathbf{B}$ denote the network admittance matrix. The generator at bus~$i\in\mathcal{G}$ has a quadratic cost function for active power generation with coefficients $c_{2i}$, $c_{1i}$, and $c_{0i}$.

Define a function for squared voltage magnitude:
\begin{equation} \label{opf_Vsq}
\left|V_i\right|^2 = f_{Vi}\left(V_d, V_q\right) := V_{di}^2 + V_{qi}^2.
\end{equation}
The power flow equations describe the network physics:
\begin{subequations}
\label{opf_balance}
\begin{align}\nonumber
P_{Gi} = f_{Pi}\left(V_d,V_q\right) := & P_{Di} + V_{di} \sum_{k=1}^n \left( \mathbf{G}_{ik} V_{dk} - \mathbf{B}_{ik} V_{qk} \right) &  &  \\[-5pt] 
\label{opf_Pbalance}  & + V_{qi} \sum_{k=1}^n \left( \mathbf{B}_{ik}V_{dk} + \mathbf{G}_{ik}V_{qk} \right), \\ \nonumber 
Q_{Gi} = f_{Qi}\left(V_d,V_q\right) := & Q_{Di} + V_{di} \sum_{k=1}^n \left( -\mathbf{B}_{ik}V_{dk} - \mathbf{G}_{ik} V_{qk}\right) \\[-5pt]
\label{opf_Qbalance} & + V_{qi} \sum_{k=1}^n \left( \mathbf{G}_{ik} V_{dk} - \mathbf{B}_{ik} V_{qk}\right).
\end{align}
\end{subequations}
Define a quadratic cost of active power generation:
\begin{equation}\label{objfunction}
f_{Ci}\left(V_d,V_q\right) := c_{2i} \left(f_{Pi}\left(V_d,V_q\right)\right)^2 + c_{1i} f_{Pi}\left(V_d,V_q\right) + c_{0i}.
\end{equation}

Each line $\left(l,m\right)\in\mathcal{L}$ is modeled by a $\Pi$ circuit with mutual admittance $y_{lm} = g_{lm} + \mathbf{j} b_{lm}$ (or, equivalently, a series impedance of $R_{lm} + \mathbf{j}X_{lm}$) and total shunt susceptance $b_{sh,lm}$. More flexible line models which include off-nominal voltage ratios and non-zero phase shifts can easily be incorporated into the proposed algorithm~\cite{illustrative_example}. Define expressions for the active, reactive, and apparent power flows on the line $\left(l,m\right) \in\mathcal{L}$:
\begin{subequations}
\begin{align}
\nonumber & f_{Plm}\left(V_d,V_q\right) := g_{lm}\left( V_{dl}^2 + V_{ql}^2\right) - g_{lm}\left(V_{dl}V_{dm} + V_{ql}V_{qm}\right)\\
\label{Plm}& \; +
b_{lm}\left(V_{dl}V_{qm} - V_{ql}V_{dm}\right), \\
\nonumber & f_{Qlm}\left(V_d,V_q\right) := -\left(b_{lm} + \frac{b_{sh,lm}}{2}\right)\left( V_{dl}^2 + V_{ql}^2\right) \\ \label{Qlm} 
& \; + b_{lm}\left(V_{dl}V_{dm} + V_{ql}V_{qm}\right) +
g_{lm}\left(V_{dl}V_{qm} - V_{ql}V_{dm}\right),\!\! \\ 
\label{Slm} & f_{Slm}\left(V_d, V_q \right) := \left(f_{Plm}\left(V_d, V_q \right)\right)^2 + \left(f_{Qlm}\left(V_d, V_q \right)\right)^2.
\end{align}
\end{subequations}

The OPF problem is
\begin{subequations}
\label{opf}
\begin{align}
\label{opf_obj} & \min_{V_d,V_q}\quad \sum_{i \in \mathcal{G}} f_{Ci}\left(V_d, V_q \right) \qquad \mathrm{subject\; to} \hspace{-160pt} & \\
\label{opf_P} &  \quad P_{Gi}^{min} \leq f_{Pi}\left(V_d,V_q\right) \leq P_{Gi}^{max} & \forall i \in \mathcal{N} \\
\label{opf_Q} &  \quad Q_{Gi}^{min} \leq f_{Qi}\left(V_d,V_q\right) \leq Q_{Gi}^{max} &  \forall i \in \mathcal{N} \\
\label{opf_V} &  \quad (V_{i}^{min})^2 \leq f_{Vi}\left(V_d,V_q\right) \leq (V_{i}^{max})^2 &  \forall i \in \mathcal{N}  \\
\label{opf_Slm} & \quad f_{Slm}\left(V_d,V_q\right) \leq \left(S_{lm}^{max}\right)^2 &  \forall \left(l,m\right) \in \mathcal{L} \\ 
\label{opf_Sml} & \quad f_{Sml}\left(V_d,V_q\right) \leq \left(S_{lm}^{max}\right)^2 & \forall \left(l,m\right) \in \mathcal{L} \\ 
\label{opf_Vref} & \quad V_{qi} = 0 & i\in\mathcal{S}
\end{align}
\end{subequations}
Constraint~\eqref{opf_Vref} sets the reference bus angle to zero.

\section{Computation of OPF Feasible Spaces}
\label{l:discretization}
Visualizing OPF feasible spaces helps researchers improve solution algorithms. To enable such visualizations, this section proposes an algorithm for provably computing the entire feasible space to within a specified discretization tolerance. The proposed algorithm discretizes certain inequality constraints to form systems of polynomial equalities, which are solved using the NPHC algorithm~\cite{guo1990,mehta2014a,SW:05,BHSW06}.

\subsection{Discretization of Inequality Constraints}
This paper discretizes certain of the OPF problem's inequality constraints to construct equality constraints in the form of power flow equations. For a set of power flow equations, load buses ($i\in\mathcal{N}\setminus\mathcal{G}$) have specified active and reactive power injections $-P_{Di}-\mathbf{j}Q_{Di}$. A single generator bus, denoted by $i\in\mathcal{S}$, is selected as the slack bus with a specified voltage phasor $V_{di} = \left|V_{i}\right|$, $V_{qi} = 0$. The active power  generation $P_{Gi}$ and squared voltage magnitudes $\left|V_{i}\right|^2$ are specified at the remaining generator buses ($i\in\mathcal{G}\setminus\mathcal{S}$). The squared voltage magnitudes at generator buses $\left|V_i\right|$, $i\in\mathcal{G}$, and active power injections at non-slack generator buses $P_{Gi}$, $i\in\mathcal{G}\setminus\mathcal{S}$, are determined using the following discretization.

Specify discretization parameters $\Delta_P$ and $\Delta_V$ for the active power injections and voltage magnitudes, respectively. The chosen discretization yields the set of power flow equations
\begin{subequations}
\vspace{-10pt}
\label{eq:pf}
\begin{align}
\label{eq:pf_Pg}
& f_{Pi}\left(V_d,V_q\right) = P_i^{min} + \eta_i \Delta_P & i\in\mathcal{G}\setminus\mathcal{S} \\
\label{eq:pf_V}
& f_{Vi}\left(V_d,V_q\right) = \left(V_{i}^{min} + \mu_i \Delta_V\right)^2 & i\in\mathcal{G}\setminus\mathcal{S}\\
& f_{Pi}\left(V_d,V_q\right) = 0 & i\in\mathcal{N}\setminus\mathcal{G}\\
& f_{Qi}\left(V_d,V_q\right) = 0 & i\in\mathcal{N}\setminus\mathcal{G}\\
\label{eq:pf_Vslack}
& V_{di} = \left(V_{i}^{min} + \mu_i \Delta_V\right) & i\in\mathcal{S} \\
& V_{qi} = 0 & i\in\mathcal{S}
\end{align}
\end{subequations}
for each combination of $\eta_i \in \left\lbrace 0,\ldots,\eta_i^{max}\right\rbrace$, $i\in\mathcal{G}\setminus\mathcal{S}$, and $\mu_i \in \left\lbrace 0,\ldots,\mu_i^{max}\right\rbrace$, $i\in\mathcal{G}$, where $\eta_i^{max} := \lfloor \frac{P_i^{max} - P_i^{min}}{\Delta_P}\rfloor$, $\mu_i^{max} := \lfloor\frac{V_i^{max} - V_i^{min}}{\Delta_V}\rfloor$, and $\lfloor\cdot\rfloor$ is the ``integer floor'' function. The number of discretization points depends on the number of generator buses $\left|\mathcal{G}\right|$; the range of the inequality constraints $P_i^{max} - P_i^{min}$, $i\in\mathcal{G}\setminus\mathcal{S}$, and $V_i^{max} - V_i^{min}$, $i\in\mathcal{G}$; and the chosen discretization parameters~$\Delta_P$ and $\Delta_V$.\footnote{To reduce the number of discretization points, the slack bus $\mathcal{S}$ is chosen as the generator bus $i\in\mathcal{G}$ with the largest value of $P_i^{max} - P_i^{min}$.}

The discretization~\eqref{eq:pf} ensues the satisfaction of the OPF problem's constraints on active power generation (other than at the slack bus) and generator voltage magnitudes as well as the load demands. Solutions to~\eqref{eq:pf} that satisfy the other inequality constraints in~\eqref{opf} are in the feasible space of the OPF problem. Thus, a ``filtering'' step is required after computing the power flow solutions at each discretization point to select the points that satisfy all the inequality constraints in~\eqref{opf}.

\subsection{Numerical Polynomial Homotopy Continuation Algorithm}
\label{l:nphc}
Ensuring the computation of the complete feasible space requires a robust algorithm for solving the power flow equations from~\eqref{eq:pf}. The Numerical Polynomial Homotopy Continuation (NPHC) algorithm~\cite{guo1990,mehta2014a,SW:05,BHSW06} is used for this purpose. 

The NPHC algorithm yields all complex solutions to systems of polynomial equations. This algorithm uses continuation to trace all the complex solutions for a ``start'' system of polynomial equations to a ``target'' system along a one-dimensional parameterization. The start system is designed such that 1) the number of complex solutions to the start system upper bounds the number of complex solutions to the target system, and 2) all solutions to the start system can be trivially computed. The NPHC algorithm guarantees that each solution to the target system is connected via a continuation trace to at least one solution of the start system~\cite{SW:05}. 

Consider a target system of $m$ quadratic equations $f_i\left(x\right) = 0$, $i=1,\ldots,m$, and variables $x\in\mathbb{C}^m$.\footnote{While NPHC is applicable to higher-order polynomials, this section focuses on quadratics to match the form of the power flow equations~\eqref{eq:pf}.} One method for constructing the start system $g\left(x\right) = 0$ uses the B\'ezout bound~\cite{SW:05} on the number of isolated complex solutions to $f\left(x\right) = 0$. The B\'ezout bound of $2^m$ suggests a start system 
\begin{align}
\label{eq:bezout_start}
& g_i\left(x\right) := a_i x_i^2 - b_i  & i=1,\ldots,m
\end{align}
where $a_i, b_i \neq 0$ are generic complex numbers. The start system $g\left(x\right) = 0$ has $2^m$ solutions of the form $x_i = \sqrt{b_i/a_i}$. Using a predictor-corrector method, the NPHC algorithm tracks all complex solutions to
\begin{equation}
\left(1-t\right)f\left(x\right) + \kappa\, t\, g\left(x\right) = 0
\end{equation}
from $t = 1$ (i.e., the start system) to $t = 0$ (i.e., the target system). The constant $\kappa$ is a randomly chosen complex number which ensures, with probability one, that the traces do not bifurcate, turn back, or cross~\cite{SW:05}.\footnote{While some traces may diverge, each solution to the target system will be reached by at least one trace beginning at a solution to the start system.} Thus, NPHC is guaranteed to find all complex solutions to the target polynomial system. 

%
%

With each bus having two constraints and two variables, $V_{di}$ and $V_{qi}$, the power flow equations~\eqref{eq:pf} for a given discretization point are a square system of polynomial equalities which can be solved with the NPHC algorithm. Only solutions with real-valued $V_d$ and $V_q$ are physically meaningful; solutions with any non-real $V_d$ or $V_q$ variables are discarded. 


The computational burden required for each solution of the NPHC algorithm depends on the number of continuation traces. When solving multiple problems that differ only in their parameter values, one approach for reducing the number of continuation traces is to compute a \emph{parameterized} homotopy. This approach solves an initial problem with generic complex parameter values (i.e., the right hand sides of~\eqref{eq:pf_Pg}, \eqref{eq:pf_V}, and \eqref{eq:pf_Vslack}). Each set of desired parameters are then solved using start systems based on the solutions to the generic set of parameters rather than~\eqref{eq:bezout_start}. Since the generic-parameter system can have significantly fewer solutions than the B\'ezout bound, fewer continuation traces are required. This effectively ``hot starts'' the NPHC algorithm for each set of parameters.

Despite the ability to speed computation for subsequent sets of parameters, the initial solution of the generic-parameter system with the B\'ezout bound can be challenging, with a requirement for $2^{2n-2}$ continuation traces. With the B\'ezout bound, NPHC is capable of solving systems with up to the 14 buses~\cite{mehta2014a}. Future work includes leveraging recently developed tighter bounds on the number of complex solutions to the power flow equations~\cite{acc2016,chen2015bounds} to speed the initial computation required for the generic-parameter system.

\section{Eliminating Infeasible Points}
\label{l:convex}
Some of the power flow equations resulting from the discretization in~\eqref{eq:pf} may be infeasible: there may not exist any real solutions or all of the real solutions may fail to satisfy the inequality constraints of~\eqref{opf}. This section proposes two screening algorithms, ``bound tightening'' and ``grid pruning'', that use Lasserre's hierarchy of convex ``moment'' relaxations~\cite{lasserre2001,pscc2014,cedric_msdp,ibm_paper} to eliminate many infeasible points. As will be described later in this section, the bound tightening algorithm improves upon the bounds on power injections, line flows, and voltage magnitudes given in the OPF problem description~\eqref{opf}. The grid pruning algorithm then eliminates infeasible points within the tightened constraints.

\subsection{Moment Relaxation Hierarchy}
The bound tightening and grid pruning algorithms employ convex relaxations to identify provably infeasible discretization points. This section describes Lasserre's moment relaxation hierarchy~\cite{lasserre2001} as applied to the OPF problem~\cite{pscc2014,ibm_opf,cedric_msdp}, with the recognition that any convex relaxation (e.g., \cite{low_tutorial,coffrin2015qc,sun2015,ghaddar2015,lavaei_tps,molzahn_holzer_lesieutre_demarco-large_scale_sdp_opf,pscc2014,cedric_msdp,molzahn_hiskens-sparse_moment_opf,ibm_opf,josz_molzahn-complex_hierarchy,bienstock2014,coffrin2016pscc}) could be used for the bound tightening and grid pruning applications to follow.



Development of the moment relaxations begins with several definitions. Define the vector of decision variables $\hat{x} \in \mathbb{R}^{2n}$:
\begin{equation}\label{xhat}
\hat{x} := \begin{bmatrix} V_{d1} & V_{d2} & \ldots & V_{dn} & V_{q1} & V_{q2} & \ldots & V_{qn} \end{bmatrix}^\intercal.
\end{equation}
A monomial is defined using an exponent vector $\alpha\in \mathbb{N}^{2n}$: $\hat{x}^{\alpha} := V_{d1}^{\alpha_1}V_{d2}^{\alpha_2}\cdots V_{qn}^{\alpha_{2n}}$. A polynomial $g\left(\hat{x}\right) := \sum_{\alpha\in \mathbb{N}^{2n}} g_{\alpha} \hat{x}^{\alpha}$, where $g_{\alpha}$ is the scalar coefficient corresponding to the monomial $\hat{x}^\alpha$.

Define a linear functional $L_{y}\left(g\right)$ which replaces the monomials $\hat{x}^{\alpha}$ in a polynomial $g\left(\hat{x}\right)$ with scalar variables $y_\alpha$:
\begin{equation}
\label{eq:Lcomp}
L_{y}\left\lbrace g \right\rbrace := \sum_{\alpha \in \mathbb{N}^{2n}} g_{\alpha} y_{\alpha}.
\end{equation}
For a matrix $g\left(\hat{x}\right)$, $L_{y}\left\lbrace g\right\rbrace$ is applied componentwise.

Consider, e.g., the vector $\hat{x} = \begin{bmatrix}V_{d1} & V_{d2} & V_{q1} & V_{q2} \end{bmatrix}^\intercal$ for a two-bus system and the polynomial $g\left(\hat{x}\right) = -\left(V_2^{min}\right)^2 + V_{d2}^2 + V_{q2}^2$. (The constraint $g\left(\hat{x}\right) \geq 0$ forces the voltage magnitude at bus~2 to be greater than or equal to $V_2^{min}$~per unit.) Then $L_{y}\left\lbrace g\right\rbrace = -\left(V_2^{min}\right)^2y_{0200} + y_{0002}$. Thus, $L_{y}\left\lbrace g \right\rbrace$ converts a polynomial $g\left(\hat{x}\right)$ to a linear function of $y$.

For the order-$\gamma$ relaxation, define a vector $x_\gamma$ consisting of all monomials of the voltages up to order $\gamma$ (i.e., $\hat{x}^{\alpha}$ such that $\left|\alpha\right| \leq \gamma$, where $\left|\,\cdot\,\right|$ is the one-norm):
\begin{align} \nonumber
x_\gamma = & \left[ \begin{array}{ccccccc} 1 & V_{d1} & \ldots & V_{qn} & V_{d1}^2 & V_{d1}V_{d2} & \ldots \end{array} \right. \\ \label{x_d}
& \qquad \left.\begin{array}{cccccc} \ldots & V_{qn}^2 & V_{d1}^3 & V_{d1}^2 V_{d2} & \ldots & V_{qn}^\gamma \end{array}\right]^\intercal
\end{align}

The relaxations are composed of positive-semidefinite-constrained \emph{moment} and \emph{localizing} matrices. The symmetric moment matrix $\mathbf{M}_{\gamma}$ has entries $y_{\alpha}$ corresponding to all monomials $\hat{x}^{\alpha}$ such that $\left|\alpha\right| \leq 2\gamma$:
\begin{equation}
\label{eq:moment}
\mathbf{M}_\gamma \left\lbrace y \right\rbrace := L_{y}\left\lbrace x_\gamma^{\vphantom{\intercal}} x_\gamma^\intercal\right\rbrace.
\end{equation}

Symmetric localizing matrices are defined for each constraint of~\eqref{opf}. For a polynomial constraint $g\left(\hat{x}\right)\geq 0$ with largest degree $\left|\alpha\right|$ among all monomials equal to $2\eta$, the localizing matrix is:
%
\begin{equation}
\label{eq:comp_local}
\mathbf{M}_{\gamma - \eta} \left\lbrace g y \right\rbrace := L_{y} \left\lbrace g\, x_{\gamma-\eta}^{\vphantom{\intercal}} x_{\gamma-\eta}^{\intercal} \right\rbrace.
\end{equation}
See~\cite{pscc2014,illustrative_example} for example moment and localizing matrices for the second-order relaxation applied to small OPF problems.

The objective functions used for the bound tightening and grid pruning algorithms in Sections~\ref{l:bt} and~\ref{l:gp} are either 1)~linear functions of the active and reactive power generation, squared voltage magnitudes, and apparent power line flows or 2)~convex quadratic functions of the active powers and squared voltage magnitudes. This section considers a general polynomial objective function $h\left(V_d,V_q\right)$ which represents a generic function in either of these forms.

The order-$\gamma$ moment relaxation is
\begin{subequations}
\label{eq:msosr}
\begin{align}
\label{eq:msosr_obj} & \min_{y} \quad L_{y}\left\lbrace h \right\rbrace \qquad \mathrm{subject\; to} \hspace{-150pt} &  \\
\label{eq:msosr_Pmin} & \mathbf{M}_{\gamma-1}\left\lbrace \left(f_{Pi} - P_i^{min}\right) y \right\rbrace \succeq 0 & \forall i\in\mathcal{N}\\
\label{eq:msosr_Pmax} &  \mathbf{M}_{\gamma-1}\left\lbrace \left(P_i^{max} - f_{Pi} \vphantom{P_i^{min}}\right) y \right\rbrace \succeq 0 & \forall i\in\mathcal{N}\\
\label{eq:msosr_Qmin} & \mathbf{M}_{\gamma-1}\left\lbrace \left(f_{Qi} - Q_i^{min}\right) y \right\rbrace \succeq 0 & \forall i\in\mathcal{N}\\
\label{eq:msosr_Qmax} &  \mathbf{M}_{\gamma-1}\left\lbrace \left(Q_i^{max} - f_{Qi}  \vphantom{P_i^{min}}\right) y \right\rbrace \succeq 0 & \forall i\in\mathcal{N}\\
\label{eq:msosr_Vmin} &  \mathbf{M}_{\gamma-1}\left\lbrace \left(f_{Vi} - \left(V_i^{min}\right)^2\right) y \right\rbrace \succeq 0 & \forall i\in\mathcal{N}\\
\label{eq:msosr_Vmax} & \mathbf{M}_{\gamma-1}\left\lbrace \left(\left(V_i^{max}\right)^2 - f_{Vi}  \vphantom{P_i^{min}}\right) y \right\rbrace \succeq 0 & \forall i\in\mathcal{N} \\
\label{eq:msosr_Slm} & \mathbf{M}_{\gamma-2}\left\lbrace \left(\left(S_{lm}^{max}\right)^2 - f_{Slm} \vphantom{P_i^{min}}\right) y \right\rbrace \succeq 0 & \forall \left(l,m\right)\in\mathcal{L} \\
\label{eq:msosr_Sml} & \mathbf{M}_{\gamma-2}\left\lbrace \left(\left(S_{ml}^{max}\right)^2 - f_{Sml} \vphantom{P_i^{min}}\right) y \right\rbrace \succeq 0 & \forall \left(l,m\right)\in\mathcal{L} \\
\label{eq:msosr_Msdp} & \mathbf{M}_\gamma \{y\} \succeq 0 & \\
\label{eq:msosr_Vref} & y_{\star\ldots\star\rho\star\ldots\star} = 0 & \rho =1,\ldots,2\gamma \\
\label{eq:msosr_y0} & y_{0\ldots 0} = 1
\end{align}
\end{subequations}
where $\rho$ in the angle reference constraint~\eqref{eq:msosr_Vref} is the index $n+k$, where $k\in\mathcal{S}$ is the index of the reference bus. Alternatively, the angle reference~\eqref{opf_Vref} can be used to eliminate all terms corresponding to $V_{qk}$, $k\in\mathcal{S}$, to reduce the problem size. Constraint~\eqref{eq:msosr_y0} corresponds to the fact that $x^0 = 1$.

For general polynomial optimization problems, the relaxation order $\gamma$ must be greater than or equal to half the largest degree of any polynomial. Objectives that are quadratic in power generation and/or squared voltage magnitudes as well as functions for apparent power line flows give rise to quartic polynomials in the voltage components, which suggests that a relaxation order $\gamma \geq 2$ is required for problems that include these functions. However, second-order cone programming (SOCP) reformulations for these functions enable the solution of~\eqref{eq:msosr} with $\gamma = 1$~\cite{lavaei_tps,illustrative_example}. Note that the first-order relaxation is equivalent to the SDP relaxation of~\cite{lavaei_tps}.

Formally, for $\gamma = 1$, the apparent power line flow limits~\eqref{opf_Slm} and~\eqref{opf_Sml} take the form of the SOCP constraints
\begin{align}
\nonumber &  S_{lm}^{max} \geq \left|\left| \begin{bmatrix} L_y\left\lbrace f_{Plm} \right\rbrace \\ L_y\left\lbrace f_{Qlm} \right\rbrace  \end{bmatrix} \right|\right|_2, & S_{lm}^{max} \geq \left|\left| \begin{bmatrix} L_y\left\lbrace f_{Pml} \right\rbrace \\ L_y\left\lbrace f_{Qml} \right\rbrace \end{bmatrix} \right|\right|_2 \hspace{-30pt}\\ \label{sdp_Slm} & & \forall \left(l,m \right) \in \mathcal{L}
\end{align}
where $\left|\left|\,\cdot\,\right|\right|_2$ is the two-norm. 
%
%
Formulation of the quartic objective function for the grid pruning algorithm is addressed in Section~\ref{l:gp}.


The relaxation~\eqref{eq:msosr} yields a single global solution if 
\begin{align}\label{rankcondition} \mathrm{rank}\big(\mathbf{M}_{\gamma}\{y\}\big) = 1.
\end{align}
The global solution $V^\ast$ is calculated using an eigendecomposition of the diagonal block of the moment matrix corresponding to the second-order monomials (i.e., $\left|\alpha\right| = 2$). Let $\sigma$ be a unit-length eigenvector corresponding to the non-zero eigenvalue $\lambda$ of $\big[\mathbf{M}_\gamma \{y\} \big]_{\left(2:2n+1,2:2n+1\right)}$. Then the globally optimal voltages are $V^\ast = \sqrt{\lambda} \sigma$. Relaxations in the moment hierarchy are guaranteed to yield the global optima of generic polynomial optimization problems at a finite relaxation order~\cite{nie2014}.

If the rank condition~\eqref{rankcondition} is satisfied, the relaxation's objective value is equal to the non-convex problem's globally optimal objective value. If the rank condition~\eqref{rankcondition} is not satisfied for some relaxation order (i.e., $\mathrm{rank}\big(\mathbf{M}_{\gamma}\{y\}\big) > 1$), the objective value of the relaxation provides a (potentially strict) lower bound on the optimal objective value for the corresponding non-convex problem. The lower bound is used in the bound tightening and grid pruning algorithms to eliminate provably infeasible points in the discretization~\eqref{eq:pf}.

\subsection{Bound Tightening}
\label{l:bt}
The bounds on the voltage magnitudes at generator buses and on the active power outputs at non-slack generator buses determine the number of discretization points in~\eqref{eq:pf}. The bounds on these quantities specified in the OPF problem may be larger than the values that are actually achievable due to the limitations imposed by other constraints. In other words, certain bounds may never be binding. It may therefore be possible to reduce the number of discretization points by determining tighter bounds on the generators' active power outputs and voltage magnitudes. This can be accomplished using a ``bound tightening'' algorithm similar to those proposed in~\cite{coffrin_tightening,sun2015,chen2015} for the purpose of determining better lower bounds on the global solutions of OPF problems.

Moment relaxations are used to tighten the OPF problem's bounds on the generators' active and reactive power outputs~\eqref{opf_P}, \eqref{opf_Q}, apparent power line flows~\eqref{opf_Slm}, \eqref{opf_Sml}, and squared voltage magnitudes~\eqref{opf_V}. Define $h_{c,\gamma}\left\lbrace f \right\rbrace$ as the solution to the following optimization problem:
\begin{align} \nonumber
& h_{c,\gamma}\left\lbrace f \right\rbrace := \max_{y} \quad L_{y}\left\lbrace cf \right\rbrace \\ \label{h} & \quad \text{subject to  } \eqref{eq:msosr_Pmin}\text{--}\eqref{eq:msosr_y0} \text{ with relaxation order } \gamma
\end{align}
where the parameter $c\in\left\lbrace -1,1\right\rbrace$ effectively determines whether the objective is to minimize or maximize, $f$ is the function corresponding to one of the OPF problem's constraints~\eqref{opf_P}--\eqref{opf_Sml}, and $\gamma$ is the specified relaxation order.

\newfloat{algorithm}{bt}{lop}
\begin{algorithm}
\caption{Bound Tightening}\label{a:bt}
\begin{algorithmic}[1]
\State \textbf{Input}: $\gamma^{max}$, upper and lower bounds $\zeta_u$ and $\zeta_{\ell}$, constraint functions $f$
\State Set $\mathcal{C}_u$ to contain all upper bound constraints
\State Set $\mathcal{C}_\ell$ to contain all lower bound constraints
\Repeat
	\ParForEach{constraint in $\mathcal{C}_u$}
		\For {$\gamma = 1,\ldots,\gamma^{max}$}
			\If {the constraint is a flow limit for line $\left(l,m\right)$}
				\If {\resizebox{170pt}{!}{$\max \left( h_{1,\gamma}\left\lbrace f_{Slm}\right\rbrace, h_{1,\gamma}\left\lbrace f_{Sml} \right\rbrace \right) < \left(S_{lm}^{max}\right)^2$}}
					\State \parbox[c][2.9em]{\dimexpr\linewidth-\algorithmicindent-80pt}{Set the flow limit for line $\left(l,m\right)$ to \\ \hspace*{5pt} $\sqrt{\max \left( h_{1,\gamma}\left\lbrace f_{Slm}\right\rbrace, h_{1,\gamma}\left\lbrace f_{Sml} \right\rbrace \right)}$}
				\EndIf
			\Else
				\If $h_{1,\gamma}\left\lbrace f \right\rbrace < \zeta_u$
					\State Update the bound: $\zeta_u \leftarrow h_{1,\gamma}\left\lbrace f \right\rbrace$
				\EndIf
			\EndIf
		\If {the rank condition~\eqref{rankcondition} is satisfied}
			\State Remove constraint from $\mathcal{C}_u$
			\Break
		\EndIf
		\EndFor
	\EndFor
	\ParForEach{constraint in $\mathcal{C}_\ell$}
		\For {$\gamma = 1,\ldots,\gamma^{max}$}
			\If $h_{\text{-}1,\gamma}\left\lbrace f \right\rbrace > \zeta_\ell$
				\State Update the bound: $\zeta_\ell \leftarrow h_{\text{-}1,\gamma}\left\lbrace f \right\rbrace$
			\EndIf
			\If {the rank condition~\eqref{rankcondition} is satisfied}
				\State Remove constraint from $\mathcal{C}_\ell$
				\Break
			\EndIf
		\EndFor
	\EndFor
\Until{no bounds are updated during this iteration}
\end{algorithmic}
\end{algorithm}

Algorithm~\ref{a:bt} describes the bound tightening approach. Given the tightest known bounds, each iteration uses~\eqref{h} to compute new bounds on the maximum and minimum achievable values of the expressions for the constrained quantities in the OPF problem~\eqref{opf}. Within an iteration, the bounds for each quantity are computed in parallel. Increasing relaxation orders of the moment hierarchy are used to determine tighter bounds. A solution to the relaxation which satisfies the rank condition~\eqref{rankcondition} yields a feasible point for the OPF problem~\eqref{opf}. No further tightening of that constraint is possible and the constraint is removed from the list of considered constraints.
%
The algorithm terminates upon reaching a fixed point where no bound tightening occurs at some iteration.

There is a subtly regarding the tightening of apparent power line flow limits. The Schur complement formulation for the the apparent power line flow limits~\eqref{sdp_Slm} cannot be maximized in an objective function. Thus, the first-order moment relaxation cannot be directly applied to tighten these limits. However, the first-order relaxation can still applied through the use of an upper bound on the apparent power line flows. Specifically, the squared line flow limits are bounded by the maximum value of the squared magnitude of the current flow multiplied by the upper bound on the squared voltage magnitude at the corresponding terminal bus. For the line $\left(l,m\right)\in\mathcal{L}$, the squared magnitude of the current flow is
\begin{align}\nonumber
& f_{Ilm}\left(V_d,V_q\right) := \left(b_{lm}^2 + g_{lm}^2\right) \left(V_{dm}^2 + V_{qm}^2\right) \\ \nonumber & \quad + b_{lm}b_{sh,lm} \left(V_{dl}V_{dm} + V_{ql} V_{qm}\right) \\ \nonumber & \quad + \left(b_{lm}^2 - b_{lm}b_{sh,lm} + b_{sh,lm}^2/4 + g_{lm}^2\right) \left(V_{dl}^2 + V_{ql}^2\right) \\ \nonumber & \quad - 2\left(b_{lm}^2+g_{lm}^2\right) \left(V_{ql}V_{qm} + V_{dl}V_{dm}\right) \\ &\quad + b_{sh,lm}g_{lm} \left(V_{dl} V_{qm} - V_{dm}V_{ql}\right)
\end{align}
The first-order relaxation can be used to obtain upper bounds on the apparent power line flow limits for the line $\left(l,m\right)\in\mathcal{L}$ by maximizing $L_y\left\lbrace\left(V_l^{max}\right)^2 f_{Ilm}\right\rbrace$ and $L_y\left\lbrace\left(V_m^{max}\right)^2 f_{Iml}\right\rbrace$. Higher-order relaxations directly formulate the expressions $L_y\left\lbrace f_{Slm} \right\rbrace$ and $L_y\left\lbrace f_{Sml} \right\rbrace$.

\subsection{Grid Pruning}
\label{l:gp}
Even the tightest possible constraints may still admit \mbox{infeasible} points in the discretization~\eqref{eq:pf}. The ``grid pruning'' algorithm described in this section often eliminates many of these infeasible points. This algorithm projects a specified point in the space of active powers and squared voltage magnitudes onto the feasible space of a convex relaxation of the OPF problem's constraints. A non-zero objective value provides the right hand side of an ellipse centered at the specified point. No feasible points for the OPF problem exist within this ellipse.

Formally, consider the optimization problem
\begin{align}\nonumber
& \phi_{\gamma}\left( P^\circ,V^\circ,\beta^\circ \right) :=  \\& \nonumber \quad \min_{y}  L_{y}\left\lbrace \sum_{i\in\mathcal{G}\setminus \mathcal{S}}  \left( f_{Pi} - P^{\circ}_i\right)^2  + \beta^\circ \sum_{i\in\mathcal{G}} \left( f_{Vi} - \left(V^{\circ}_i\right)^2\right)^2 \right\rbrace \\ & \label{proj} \qquad \text{subject to  } \eqref{eq:msosr_Pmin}\text{--}\eqref{eq:msosr_y0} \text{ with relaxation order }\gamma
\end{align}
where $P^\circ\in\mathbb{R}^n$ and $V^\circ\in\mathbb{R}^n$ are vectors of parameters specifying a point in the space of active powers and voltage magnitudes, the parameter $\beta^\circ$ specifies a scalar coefficient that weights distances in active power to distances in squared voltage magnitude, and $\gamma$ is the relaxation order.

For $\phi_{1}\left( P^\circ,V^\circ,\beta^\circ \right)$, the objective in~\eqref{proj} minimizes an auxiliary variable $\omega$ with the SOCP constraints
\begin{align}
\label{proj_socp} & \omega \geq \left|\left| \begin{bmatrix} \left(L_y\left\lbrace F_{Pi} \right\rbrace\right)^\intercal & \left(L_y\left\lbrace F_{Vi} \right\rbrace\right)^\intercal \end{bmatrix}^\intercal \right|\right|_2
\end{align}
where $F_{Pi}$ and $F_{Vi}$ are the vectors containing $f_{Pi}$, $\forall i\in\mathcal{G}\setminus\mathcal{S}$, and $f_{Vi}$, $\forall i\in\mathcal{G}$, respectively.

A solution to~\eqref{proj} with $\phi_{\gamma}\left( P^\circ,V^\circ,\beta^\circ \right) > 0$ provides the right hand side of an ellipse in the space of active powers $P$ and voltage magnitudes $V$ that is centered at $P^\circ$ and $V^\circ$ with the weighting between squared voltages and active power generation described by $\beta^\circ$:
\begin{equation}\label{ellipse}
\sum_{i\in\mathcal{G}\setminus \mathcal{S}}  \left( P_i - P^{\circ}_i\right)^2  + \beta^\circ \sum_{i\in\mathcal{G}} \left( \left(V_i\right)^2 - \left(V^{\circ}_i\right)^2\right)^2\! <\! \phi_{\gamma}\left( P^\circ,V^\circ,\beta^\circ \right)
\end{equation}
Points satisfying~\eqref{ellipse} are infeasible for the OPF problem~\eqref{opf}.\footnote{If a point satisfying~\eqref{ellipse} were feasible for the OPF problem~\eqref{opf}, it would be included in the feasible space of the moment relaxation~\eqref{proj}, resulting in an objective value less than $\phi_{\gamma}\left( P^\circ,V^\circ,\beta^\circ \right)$.}

The grid pruning method in Algorithm~\ref{a:gp} uses~\eqref{ellipse} to eliminate infeasible discretization points. Consider two discretizations of the form~\eqref{eq:pf}: a ``dense'' discretization with parameters $\hat{\Delta}_P$ and $\hat{\Delta}_V$, which is denoted by $\mathcal{D}_d$ with points $\hat{P}^\circ$ and $\hat{V}^\circ$, and a ``sparse'' discretization with parameters $\overline{\Delta}_P > \hat{\Delta}_P$ and $\overline{\Delta}_V > \hat{\Delta}_V$, which is denoted by $\mathcal{D}_s$ with points $\overline{P}^\circ$ and $\overline{V}^\circ$. The dense discretization represents the feasible space of the OPF problem while the sparse discretization provides the specified points in the grid pruning algorithm. Algorithm~\ref{a:gp} solves~\eqref{proj} at each point in sparse discretization. For each solution with $\phi_{\gamma}\left( \overline{P}^\circ,\overline{V}^\circ,\beta^\circ \right) > 0$, all points $\mathcal{D}_d$ which satisfy~\eqref{ellipse} are infeasible and therefore eliminated.

This process is repeated for various values of the weighting parameter $\beta^\circ$. The choice of different weighting parameters changes the shape of the ellipse~\eqref{ellipse} and can therefore result in the elimination of additional infeasible points.


Any solution to~\eqref{proj} which satisfies the rank condition~\eqref{rankcondition} is feasible for the OPF problem~\eqref{opf}. Higher-order relaxations are not required for any point in the sparse discretization that yields a solution satisfying~\eqref{rankcondition}.

\newfloat{algorithm}{bt}{lop}
\begin{algorithm}
\caption{Grid Pruning}\label{a:gp}
\begin{algorithmic}[1]
\State \textbf{Input}: scalar $\gamma^{max}$, vector $\beta$, dense discretization $\mathcal{D}_d$ defined using~\eqref{eq:pf} with $\hat{\Delta}_P$ and $\hat{\Delta}_V$ yielding points $\hat{P}^\circ$ and $\hat{V}^\circ$, sparse discretization $\mathcal{D}_s$ defined using~\eqref{eq:pf} with $\overline{\Delta}_P > \hat{\Delta}_P$ and $\overline{\Delta}_V > \hat{\Delta}_V$ yielding points $\overline{P}^\circ$ and $\overline{V}^\circ$
\ForEach {$\beta^\circ\in\beta$}
	\State Set $\mathcal{D}_{s,\beta^\circ} \leftarrow \mathcal{D}_s$
	\For {$\gamma = 1,\ldots,\gamma^{max}$}	
		\ParForEach {point in $\mathcal{D}_{s,\beta^\circ}$}
			\State Compute $\phi_{\gamma}\left( \overline{P}^\circ,\overline{V}^\circ,\beta^\circ \right)$ with~\eqref{proj}
			\State \parbox[c][4.4em]{\dimexpr\linewidth-\algorithmicindent-30pt}{Eliminate points in  $\mathcal{D}_d$ satisfying~\eqref{ellipse} with \\ \hspace*{5pt} right hand side $\phi_{\gamma}\left( \overline{P}^\circ,\overline{V}^\circ,\beta^\circ \right)$, $P := \hat{P}^\circ$, \\ \hspace*{5pt} $V := \hat{V}^\circ$, $P^\circ := \overline{P}^\circ$, and $V^\circ := \overline{V}^\circ$}
			\If {the rank condition~\eqref{rankcondition} is satisfied}
				\State Remove this point from $\mathcal{D}_{s,\beta^\circ}$
			\EndIf
		\EndFor
	\EndFor
\EndFor
\end{algorithmic}
\end{algorithm}

\subsection{Feasible Space Computation Algorithm}
Algorithm~\ref{a:fs} describes the method for computing an OPF feasible space. First, Algorithm~\ref{a:bt} tightens the constraint bounds and then Algorithm~\ref{a:gp} eliminates provably infeasible points within the tighter bounds. The NPHC algorithm is applied (in parallel) to solve the power flow equations corresponding to the remaining discretization points. Finally, the resulting real power flow solutions are filtered to select only those satisfying all the constraints in the OPF problem~\eqref{opf}.

\newfloat{algorithm}{bt}{lop}
\begin{algorithm}
\caption{OPF Feasible Space Computation}\label{a:fs}
\begin{algorithmic}[1]
\State \textbf{Input}: OPF constraint bounds and functions, scalar $\gamma^{max}$, vector $\beta$, dense discretization parameters $\hat{\Delta}_P$ and $\hat{\Delta}_V$, sparse discretization parameters $\overline{\Delta}_P$ and $\overline{\Delta}_V$
\State Tighten bounds using Algorithm~\ref{a:bt} with relaxations up to order $\gamma^{max}$
\State Save any resulting solutions that satisfy~\eqref{rankcondition}
\State Construct dense and sparse discretizations, $\mathcal{D}_d$ and $\mathcal{D}_s$, using~\eqref{eq:pf} with $\hat{\Delta}_P$, $\hat{\Delta}_V$ and $\overline{\Delta}_P$, $\overline{\Delta}_V$, respectively
\State Prune $\mathcal{D}_d$ using Algorithm~\ref{a:gp} with $\gamma^{max}$, $\beta$, $\mathcal{D}_s$, and $\mathcal{D}_d$
\State Save any resulting solutions that satisfy~\eqref{rankcondition}
\ParForEach {discretization point in $\mathcal{D}_d$}
	\State Solve the power flow equations~\eqref{eq:pf} using NPHC
\EndFor
\State Filter the power flow solutions satisfying all constraints~\eqref{opf}
\State \textbf{Output:} Filtered power flow solutions augmented with the rank-one solutions obtained from Algorithms~\ref{a:bt} and~\ref{a:gp}
\end{algorithmic}
\end{algorithm}

\section{Example Test Cases}
\label{l:examples}
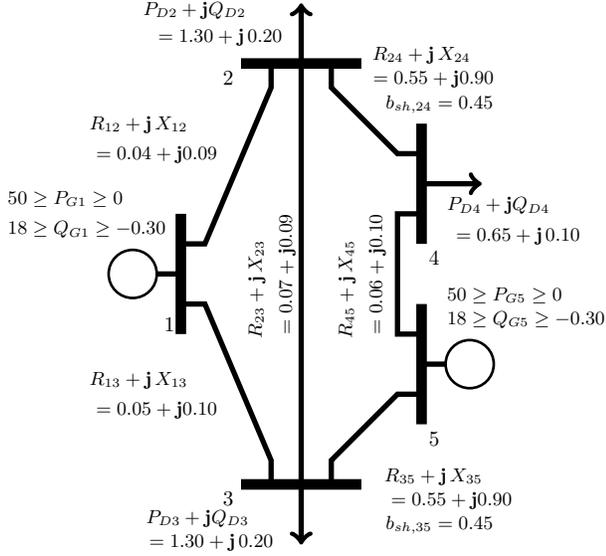
\begin{figure}[tb]
\centering
\begin{circuitikz}[scale=0.8, transform shape]
\ctikzset{bipoles/length=0.6cm}

\path[draw,line width=4pt] (0,-1) -- (0,1);
\draw (-0.4,-0.6) node[below right] {$1$};
\path[draw,line width=2pt] (-0.4,0) -- (0,0);
\draw[line width=1] (-0.8,0) circle (0.4);
\draw (-3,1.5) node[below right] {\small $50 \geq P_{G1} \geq 0$};
\draw (-3,1.0) node[below right] {\small $18\geq Q_{G1} \geq -0.30$};

\path[draw,line width=4pt] (1,3.5) -- (3,3.5);
\draw (1,3.5) node[below left] {$2$};
\path[draw,->,line width=2pt] (2,3.5) -- (2,4.5);
\draw (-0.75,4.65) node[below right] {\small $P_{D2} + \mathbf{j} Q_{D2}$};
\draw (-0.5,4.2) node[below right] {\small $=1.30 + \mathbf{j}\, 0.20$};

\path[draw,line width=4pt] (4,0.5) -- (4,2.5);
\draw (4,0.5) node[below right] {$4$};
\path[draw,->,line width=2pt] (4,1.5) -- (5,1.5);
\draw (4.3,1.4) node[below right] {\small $P_{D4} + \mathbf{j} Q_{D4}$};
\draw (4.45,0.9) node[below right] {\small $=0.65 + \mathbf{j}\, 0.10$};

\path[draw,line width=4pt] (4,-0.5) -- (4,-2.5);
\draw (4,-2.5) node[below right] {$5$};
\path[draw,line width=2pt] (4,-1.5) -- (4.4,-1.5);
\draw[line width=1] (4.8,-1.5) circle (0.4);
\draw (4.3,-0.1) node[below right] {\small $50 \geq P_{G5} \geq 0$};
\draw (4.3,-0.5) node[below right] {\small $18 \geq Q_{G5} \geq -0.30$};

\path[draw,line width=4pt] (1,-3.5) -- (3,-3.5);
\draw (1,-3.5) node[below left] {$3$};
\path[draw,->,line width=2pt] (2,-3.5) -- (2,-4.5);
\draw (-0.7,-3.8) node[below right] {\small $P_{D3} + \mathbf{j} Q_{D3}$};
\draw (-0.65,-4.2) node[below right] {\small $=1.30 + \mathbf{j}\, 0.20$};

\path[draw,line width=2pt] (0,0.5) -- (0.4,0.5) -- (1.5,3.1) -- (1.5,3.5);
\draw (-0.7,2.75) node[anchor=north] {\small $R_{12}+\mathbf{j}\,X_{12}$};
\draw (-0.4,2.25) node[anchor=north] {\small $=0.04+\mathbf{j}0.09$};

\path[draw,line width=2pt] (0,-0.5) -- (0.4,-0.5) -- (1.5,-3.1) -- (1.5,-3.5);
\draw (-0.7,-1.5) node[anchor=north] {\small $R_{13}+\mathbf{j}\,X_{13}$};
\draw (-0.4,-2.0) node[anchor=north] {\small $=0.05+\mathbf{j}0.10$};

\path[draw,line width=2pt] (2,3.5) -- (2,-3.5);
\draw (1,-0.25) node[rotate=90,anchor=north] {\small $R_{23}+\mathbf{j}\,X_{23}$};
\draw (1.5,-0.) node[rotate=90,anchor=north] {\small $=0.07+\mathbf{j}0.09$};

\path[draw,line width=2pt] (2.5,3.5) -- (2.5,3.1) -- (3.6,2) -- (4,2);
\draw (4,3.9) node[rotate=0,anchor=north] {\small $R_{24}+\mathbf{j}\,X_{24}$};
\draw (4.2,3.5) node[rotate=0,anchor=north] {\small $=0.55+\mathbf{j}0.90$};
\draw (4.3,3.1) node[rotate=0,anchor=north] {\small $b_{sh,24}=0.45$};

\path[draw,line width=2pt] (4,1) -- (3.6,1) -- (3.6,-1) -- (4,-1);
\draw (2.5,-0.25) node[rotate=90,anchor=north] {\small $R_{45}+\mathbf{j}\,X_{45}$};
\draw (3,-0.) node[rotate=90,anchor=north] {\small $=0.06+\mathbf{j}0.10$};

\path[draw,line width=2pt] (4,-2) -- (3.6,-2) -- (2.5,-3.1) -- (2.5,-3.5);
\draw (4.2,-3.1) node[rotate=0,anchor=north] {\small $R_{35}+\mathbf{j}\,X_{35}$};
\draw (4.5,-3.5) node[rotate=0,anchor=north] {\small $=0.55+\mathbf{j}0.90$};
\draw (4.3,-3.9) node[rotate=0,anchor=north] {\small $b_{sh,35}=0.45$};

\end{circuitikz}
\caption{Five-bus system from~\cite{bukhsh_tps} with impedances and powers in per unit}
\label{f:WB5}
\end{figure}

\begin{figure}[ht]
\centering
\begin{circuitikz}[scale=0.8, transform shape]
\ctikzset{bipoles/length=0.6cm}



\path[draw,line width=4pt] (-5,1.25) -- (-3,1.25);
\draw (-5,1.24) node[below right] {$1$};
\path[draw,line width=2pt] (-4.25,1.25) -- (-4.25,1.25-0.4);
\draw[line width=1] (-4.25,1.25-0.8) circle (0.4);
\draw (-3.75,0.8) node[below right] {\small $2.50 \geq P_{G1} \geq 0$};
\draw (-3.75,0.5) node[below right] {\small $3.00 \geq Q_{G5} \geq -0.30$};

\path[draw,line width=2pt] (0,2) -- (0,1) -- (-3.5,1) -- (-3.5,1.25);
\draw (0.2,1.4) node[below right] {\small $R_{14}+\mathbf{j}\,X_{14}$};
\draw (0.2,1.) node[below right] {\small $=0 + \mathbf{j}\, 0.0576$};

\path[draw,line width=4pt] (-1,2) -- (1,2);
\draw (-1,1.9) node[below right] {$4$};

\path[draw,line width=4pt] (-3,4) -- (-1,4);
\draw (-1,4) node[below left] {$9$};
\path[draw,->,line width=2pt] (-2.5,4) -- (-2.5,3);
\draw (-5,3.7) node[below right] {\small $P_{D9} + \mathbf{j}Q_{D9}$};
\draw (-5,3.2) node[below right] {\small $ = 0.75 + \mathbf{j}0.30$};

\path[draw,line width=4pt] (1,4) -- (3,4);
\draw (1,4) node[below right] {$5$};
\path[draw,->,line width=2pt] (2.5,4) -- (2.5,3);
\draw (2.7,3.7) node[below right] {\small $P_{D5} + \mathbf{j}Q_{D5}$};
\draw (2.7,3.2) node[below right] {\small $ = 0.54 + \mathbf{j}0.18$};

\path[draw,line width=4pt] (-3,6) -- (-1,6);
\draw (-1,6) node[below left] {$8$};

\path[draw,line width=4pt] (1,6) -- (3,6);
\draw (1,6) node[below right] {$6$};

\path[draw,line width=4pt] (-1,8.5) -- (1,8.5);
\draw (1,8.4) node[below left] {$7$};
\draw (0.3,9.7) node[below right] {\small $P_{D7} + \mathbf{j}Q_{D7}$};
\draw (0.3,9.3) node[below right] {\small $ = 0.60 + \mathbf{j}0.21$};
\path[draw,->,line width=2pt] (0,8.5) -- (0,9.5);

\path[draw,line width=4pt] (-5,8) -- (-3,8);
\draw (-5,7.9) node[below right] {$2$};
\path[draw,line width=2pt] (-4,8) -- (-4,8.4);
\draw[line width=1] (-4,8.8) circle (0.4);
\draw (-5.5,10.3) node[below right] {\small $3.00 \geq P_{G2} \geq 0.10$};
\draw (-5.5,9.8) node[below right] {\small $3.00 \geq Q_{G2} \geq -0.05$};

\path[draw,line width=4pt] (5,8) -- (3,8);
\draw (5,7.9) node[below left] {$3$};
\path[draw,line width=2pt] (4,8) -- (4,8.4);
\draw[line width=1] (4,8.8) circle (0.4);
\draw (2.5,10.3) node[below right] {\small $2.70 \geq P_{G2} \geq 0.10$};
\draw (2.5,9.8) node[below right] {\small $3.00 \geq Q_{G2} \geq -0.05$};

\path[draw,line width=2pt] (-0.5,2) -- (-0.5,2.3) -- (-1.5,3.7) -- (-1.5,4);
\draw (-3.5,2.8) node[below right] {\small $R_{49}+\mathbf{j}\,X_{49}$};
\draw (-4.0,2.4) node[below right] {\small $=0.0170 + \mathbf{j}\, 0.0920$};
\draw (-3.6,2) node[below right] {\small $b_{sh,49}=0.1580$};

\path[draw,line width=2pt] (0.5,2) -- (0.5,2.3) -- (1.5,3.7) -- (1.5,4);
\draw (1.2,2.8) node[below right] {\small $R_{78}+\mathbf{j}\,X_{78}$};
\draw (1.3,2.4) node[below right] {\small $=0.0170 + \mathbf{j}\, 0.0920$};
\draw (1.6,2) node[below right] {\small $b_{sh,78}=0.1580$};

\path[draw,line width=2pt] (-2,4) -- (-2,6);
\draw (-5,5.5) node[below right] {\small $R_{89}+\mathbf{j}\,X_{89}$};
\draw (-5,5.1) node[below right] {\small $=0.0320 + \mathbf{j}\, 0.1610$};
\draw (-5,4.7) node[below right] {\small $b_{sh,89}=0.3060$};

\path[draw,line width=2pt] (2,4) -- (2,6);
\draw (2.2,5.5) node[below right] {\small $R_{56}+\mathbf{j}\,X_{56}$};
\draw (2.2,5.1) node[below right] {\small $=0.0390 + \mathbf{j}\, 0.1700$};
\draw (2.2,4.7) node[below right] {\small $b_{sh,56}=0.3580$};

\path[draw,line width=2pt] (-1.5,6) -- (-1.5,6.3) -- (-0.5,7.7) -- (-0.5,8.5);
\draw (-3,8) node[below right] {\small $R_{78}+\mathbf{j}\,X_{78}$};
\draw (-3.7,7.6) node[below right] {\small $=0.0085 + \mathbf{j}\, 0.0720$};
\draw (-3.5,7.1) node[below right] {\small $b_{sh,78}=0.1490$};

\path[draw,line width=2pt] (1.5,6) -- (1.5,6.3) -- (0.5,7.7) -- (0.5,8.5);
\draw (1.,8) node[below right] {\small $R_{67}+\mathbf{j}\,X_{67}$};
\draw (1.1,7.6) node[below right] {\small $=0.119 + \mathbf{j}\, 0.1008$};
\draw (1.2,7.1) node[below right] {\small $b_{sh,67}=0.2090$};

\path[draw,line width=2pt] (-2.5,6) -- (-2.5,6.5) -- (-4,6.5) -- (-4,8);
\draw (-5.5+0.5,6.5-0.15) node[below right] {\small $R_{28}+\mathbf{j}\,X_{28}$};
\draw (-5.5+0.5,6.1-0.15) node[below right] {\small $=0 + \mathbf{j}\, 0.0625$};

\path[draw,line width=2pt] (2.5,6) -- (2.5,6.5) -- (4,6.5) -- (4,8);
\draw (4-0.75,6.5-0.15) node[below right] {\small $R_{36}+\mathbf{j}\,X_{36}$};
\draw (4-0.75,6.1-0.15) node[below right] {\small $=0 + \mathbf{j}\, 0.0586$};

\end{circuitikz}
\caption{Nine-bus system from~\cite{bukhsh_tps} with impedances and powers in per unit}
\label{f:case9mod}
\vspace*{-1em}
\end{figure}
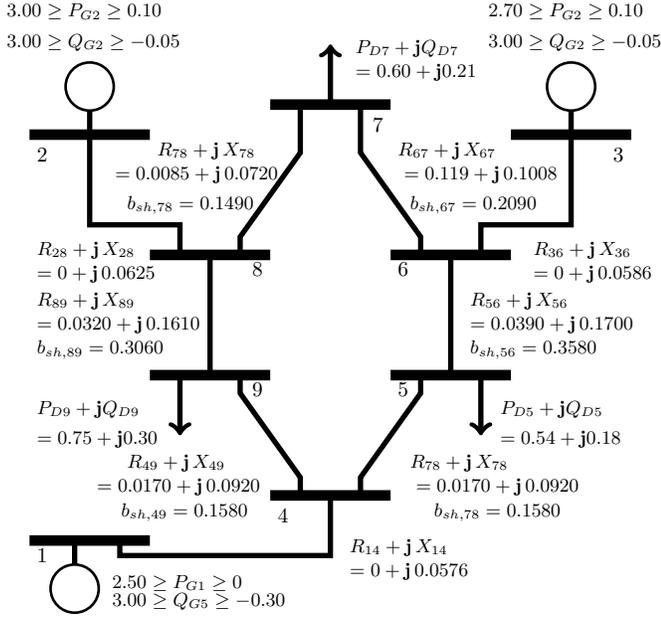

This section applies Algorithm~\ref{a:fs} to two small OPF test cases which have multiple local optima~\cite{bukhsh_tps}. The five-bus system ``WB5'' has the one-line diagram in Fig.~\ref{f:WB5}. The voltage magnitudes in WB5 are constrained to the range $\left| V_i\right| \in \left[0.95,1.05\right]$~per~unit and there are no line flow limits. The nine-bus system ``case9mod'' has the one-line diagram in Fig.~\ref{f:case9mod}. The voltage magnitudes in case9mod are constrained to the range $\left| V_i\right| \in \left[0.90,1.10\right]$~per~unit and limits on the apparent power line flows are $250$~MVA for all lines except for $\left(5,6\right)$ and $\left(6,7\right)$, which are limited to $150$~MVA, and $\left(3,6\right)$, which is limited to $300$~MVA. Both test cases use a $100$~MVA base.

Both WB5 and case9mod challenge a variety of optimization algorithms. Local solvers with a variety of reasonable initializations often converge to suboptimal local solutions in these problems. The SDP relaxation of~\cite{lavaei_tps} is not exact for either test case. Conversely, the second-order moment relaxation finds the global solution to both problems.

Algorithm~\ref{a:fs} is run for each of these systems using $\gamma^{max}$ = 2 for the bound tightening and grid pruning algorithms. The implementation uses \mbox{MATLAB} with YALMIP \mbox{2015.06.26}~\cite{yalmip} and BertiniLab~v.1.5~\cite{bertinilab}, the SDP solver in Mosek \mbox{7.1.0.28}, and Bertini~v1.4.1~\cite{BHSW06}. The Fusion cluster at Argonne National Laboratory was used for the NPHC computations.

\begin{figure}[tb]
\centering
\includegraphics[totalheight=0.257\textheight]{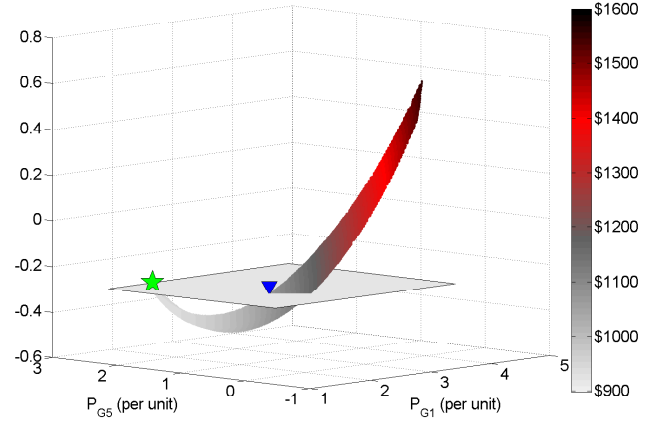}
\caption{Feasible Space for the five-bus system from~\cite{bukhsh_tps}. The colors represent the generation cost. The gray plane shows the lower reactive power limit $Q_{G5} \geq -0.30$~per~unit. This limit splits the feasible space into the two disconnected components which are above the gray plane. The green star shows the global solution and the blue triangle indicates a local optimum.}
\label{f:WB5_fs}
\end{figure}

Fig.~\ref{f:WB5_fs} shows a projection of the feasible space for WB5 in terms of the active power generations $P_{G1}$ and $P_{G5}$ and the reactive power generation $Q_{G5}$. The colors represent the generation cost corresponding to the specified objective function, $400 P_{G1} + 100 P_{G5}$. The lower limit $Q_{G5} \geq -0.30$~per~unit is shown by the gray plane. The feasible space is composed of the two disconnected components that lie above this plane. The global solution is shown by the green star at $\left(P_{G1},P_{G5},Q_{G5}\right) = \left(1.81,\, 2.21,\,-0.30 \right)$~per~unit. The blue triangle at $\left(P_{G1},P_{G5},Q_{G5}\right) = \left(2.46,\, 0.98,\,-0.30\right)$~per~unit denotes a local solution with an objective value that is $14.34$\% greater than that of the global solution.

\begin{figure}[tb]
\centering
\includegraphics[totalheight=0.26\textheight]{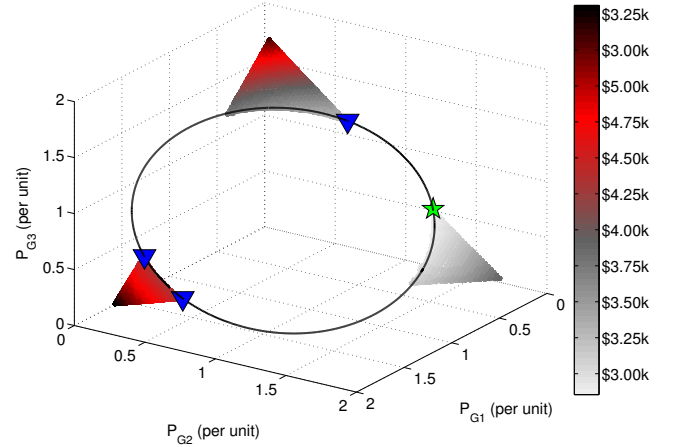}
\caption{Feasible Space for the nine-bus system from~\cite{bukhsh_tps}. The colors represent the generation cost. The feasible space is split into three disconnected components by the black line, which signifies the set of points for which the limits $Q_{G1} \geq -0.05$, $Q_{G2} \geq -0.05$, $Q_{G3} \geq -0.05$, and $\left|V_{9}\right| \geq 0.9$~per~unit are simultaneously binding. The green star shows the global solution and the blue triangles indicate local optima.}
\label{f:case9mod_fs}
\end{figure}

The feasible space for WB5 shown in Fig.~\ref{f:WB5_fs} was constructed with discretization parameters $\hat{\Delta}_P = 1$~MW and $\hat{\Delta}_V = 0.001$~per~unit. Bound tightening (Algorithm~\ref{a:bt}) eliminated 98.65\% of the points resulting from the original OPF problem's bounds. Grid pruning (Algorithm~\ref{a:gp}) using a sparse discretization with parameters $\overline{\Delta}_P = 5$~MW, $\overline{\Delta}_V = 0.005$~per~unit, and $\beta = 1$ removed 76.46\% of the remaining points in the bound-tightened problem. Thus, the bound tightening and grid pruning algorithms removed a total of 99.68\% of the initial discretization points, which suggests that the bounds specified for the OPF problem poorly represent the actual feasible space. After applying the bound tightening and grid pruning algorithms, a total of $1.62\times 10^5$ points were solved with the NPHC algorithm, with 76.46\% of these points satisfying the OPF constraints and therefore included in the feasible space. Initial solution of the parameterized NPHC algorithm described in Section~\ref{l:nphc} required 3 seconds. Each subsequent NPHC solve required approximately 0.5 seconds.

\begin{table}[tb]
\centering
\caption{Generation cost functions for the nine-bus system from~\cite{bukhsh_tps}}
\label{t:case9mod_cost}
\begin{tabular}{|c|c|c|c|}
\hline 
\textbf{Bus} & $\mathbf{c_{2i}}$ [\textbf{\$}/$\mathbf{\text{\textbf{(per~unit-hr)}}^2}$] & $\mathbf{c_{1i}}$ [\textbf{\$/(per~unit-hr)}] & $\mathbf{c_{0i}}$ [\textbf{\$/hr}]  \\ \hline\hline
1 & 1100  & 500 & 150 \\
2 & 85    & 120 & 600 \\
3 & 122.5 & 100 & 335 \\ \hline
\end{tabular}
\end{table}

Fig.~\ref{f:case9mod_fs} shows a projection of the feasible space for case9mod in terms of the active power generations $P_{G1}$, $P_{G2}$, and $P_{G3}$. The colors represent the generation cost in terms of the specified objective function, which has coefficients given in Table~\ref{t:case9mod_cost}. The feasible space has three disconnected components. 
%
The green star at $\left(P_{G1},P_{G2},P_{G3}\right) = \left(0.10,\,1.254,\,0.570\right)$~per~unit shows the global solution. The blue triangles at $\left(P_{G1},P_{G2},P_{G3}\right) = \left(0.10,\,0.648,\,1.178\right)$, $\left(1.432,\,0.378,\,0.10\right)$, and $\left(1.422,\,0.10,\,0.388\right)$~per~unit denote the three local optima, which have objective values that are $10.05$\%, $37.52$\%, and $38.13$\%, respectively, greater than the that of the global optimum.

The feasible space shown in Fig.~\ref{f:case9mod_fs} is cut by the ellipse denoted by the black line. This ellipse is comprised of points for which the lower voltage magnitude constraint at bus~9 as well as the lower reactive power limits on the generators are all binding (i.e., $\left|V_9\right|^2 = \left(0.90\right)^2$, $Q_{G1} = Q_{G2} = Q_{G3} = -0.05$~per~unit).\footnote{The points which also satisfy the other constraints in~\eqref{opf} are included in the feasible space, whereas the remainder of the black line is infeasible.} In other words, the lower voltage magnitude and lower reactive power generation limits interact to yield a disconnected feasible space. The points in the three different ``corner'' regions of Fig.~\ref{f:case9mod_fs} correspond to generator outputs that are very different active power but similar in reactive power.

Constructing the feasible space for case9mod started with a relatively sparse discretization of $\hat{\Delta}_P = 10$~MW and $\hat{\Delta}_V = 0.005$~per~unit to identify the three disconnected components of the feasible space. This facilitated multiple computations with Algorithm~\ref{a:fs} for adjoining subregions of the feasible space, with each subregion containing one of the three disconnected components. The bound tightening performed by Algorithm~\ref{a:bt} was significantly more effective when applied to each subregion, which enabled computation with a denser discretization of $\hat{\Delta}_P = 2$~MW and $\hat{\Delta}_V = 0.002$~per~unit. To improve fidelity near the dashed ellipse in Fig.~\ref{f:case9mod_fs}, a variety of smaller regions were considered with discretization tolerances up to $\hat{\Delta}_P = 1$~MW and $\hat{\Delta}_V = 0.0005$~per~unit. The grid pruning algorithm used $\beta = \left\lbrace 100, 10, 1\right\rbrace$ and a sparse discretization with parameters $\overline{\Delta}_P  = 20$~MW and $\overline{\Delta}_V  = 0.02$~per~unit. Overall, the bound tightening algorithm eliminated 99.96\% of the initial discretization points. The grid pruning algorithm removed 96.77\% of the remaining points. Thus, the bound tightening and grid pruning algorithms removed a total of 99.9987\% of the initial points in the discretization. After applying the bound tightening and grid pruning algorithms, there were $1.74\times 10^6$ remaining points which were solved (in parallel) with NPHC. Of these, 2.55\% of the NPHC solutions were feasible (i.e., passed the filtering in the last step of Algorithm~\ref{a:fs}), which suggests that it may be possible to further improve the detection of infeasible points. The initial parameterized NPHC solution required 740 seconds. Each subsequent NPHC solve required approximately 1.4 seconds.

Observe that the system parameters in Figs.~\ref{f:WB5} and~\ref{f:case9mod} are reasonable (e.g., all lines have resistance-to-reactance ratios less than one, all loads have power factors greater than $0.9$, the voltage magnitudes are constrained to be near their nominal values). Despite this, the feasible spaces for the corresponding problems exhibit significant non-convexity.

These examples illustrate that the challenges associated with certain OPF problems are strongly related to the voltage magnitude and reactive power limits. For WB5 and case9mod, binding reactive power constraints result in disconnected feasible spaces. The disconnected components contain local optima that are significantly inferior to the global optima. Disconnected feasible spaces may also result from binding apparent power line flow limits~\cite{hicss2014}. Further characterizing the physical characteristics which give rise to challenging feasible spaces is an important future research direction that will be informed by the proposed algorithm.

\section{Conclusion and Future Work}
\label{l:conclusion}
This paper has proposed an algorithm for computing the feasible spaces of small OPF problems. This algorithm discretizes certain of the OPF problem's inequality constraints to a set of power flow equations. The Numerical Polynomial Homotopy Continuation (NPHC) algorithm is used to reliably solve the power flow equations at each discretization point. The power flow solutions which satisfy all OPF constraints are included in the feasible space. Thus, the proposed algorithm is guaranteed to compute the entire feasible space to within a specified discretization tolerance. Bound tightening and grid pruning algorithms improve computational tractability by using convex moment relaxations to eliminate infeasible points. 

Future work includes computational improvements, such as integration with the software Paramotopy~\cite{paramotopy} to reduce overhead and exploitation of network structure~\cite{acc2016,chen2015bounds} to reduce the initial solution time for the parameterized NPHC algorithm. Future work also includes applying the proposed algorithm to other test cases to further characterize the physical features that are associated with challenging OPF problems.

\section*{Acknowledgment}
Discussions with Drs.~D.~Mehta, M.~Niemerg, and T.~Chen regarding NPHC are gratefully acknowledged, as is the advice from Ms.~N.-J. Simon regarding the display of the figures. Use of Fusion, a high-performance computing cluster operated by the Laboratory Computing Resource Center at Argonne National Laboratory, is also appreciated.




\bibliographystyle{IEEEtran}
\bibliography{IEEEabrv,vis_opf}

\begin{thebibliography}{10}
\providecommand{\url}[1]{#1}
\csname url@samestyle\endcsname
\providecommand{\newblock}{\relax}
\providecommand{\bibinfo}[2]{#2}
\providecommand{\BIBentrySTDinterwordspacing}{\spaceskip=0pt\relax}
\providecommand{\BIBentryALTinterwordstretchfactor}{4}
\providecommand{\BIBentryALTinterwordspacing}{\spaceskip=\fontdimen2\font plus
\BIBentryALTinterwordstretchfactor\fontdimen3\font minus
  \fontdimen4\font\relax}
\providecommand{\BIBforeignlanguage}[2]{{%
\expandafter\ifx\csname l@#1\endcsname\relax
\typeout{** WARNING: IEEEtran.bst: No hyphenation pattern has been}%
\typeout{** loaded for the language `#1'. Using the pattern for}%
\typeout{** the default language instead.}%
\else
\language=\csname l@#1\endcsname
\fi
#2}}
\providecommand{\BIBdecl}{\relax}
\BIBdecl

\bibitem{bukhsh_tps}
W.~Bukhsh, A.~Grothey, K.~McKinnon, and P.~Trodden, ``{Local Solutions of the
  Optimal Power Flow Problem},'' \emph{IEEE Trans. Power Syst.}, vol.~28,
  no.~4, pp. 4780--4788, 2013.

\bibitem{lavaei_tps}
J.~Lavaei and S.~Low, ``{Zero Duality Gap in Optimal Power Flow Problem},''
  \emph{IEEE Trans. Power Syst.}, vol.~27, no.~1, pp. 92--107, 2012.

\bibitem{bienstock2015nphard}
D.~Bienstock and A.~Verma, ``{Strong NP-hardness of AC Power Flows
  Feasibility},'' \emph{arXiv:1512.07315}, Dec. 2015.

\bibitem{NPhard}
K.~Lehmann, A.~Grastien, and P.~{Van Hentenryck}, ``{AC-Feasibility on Tree
  Networks is NP-Hard},'' \emph{IEEE Trans. Power Syst.}, vol.~31, no.~1, pp.
  798--801, Jan. 2016.

\bibitem{carpentier}
J.~Carpentier, ``{Contribution to the Economic Dispatch Problem},'' \emph{Bull.
  Soc. Franc. Elect}, vol.~8, no.~3, pp. 431--447, 1962.

\bibitem{opf_litreview1993IandII}
J.~Momoh, R.~Adapa, and M.~El-Hawary, ``{A Review of Selected Optimal Power
  Flow Literature to 1993. Parts I and II},'' \emph{IEEE Trans. Power Syst.},
  vol.~14, no.~1, pp. 96--111, Feb. 1999.

\bibitem{ferc4}
A.~Castillo and R.~O'Neill, ``{Survey of Approaches to Solving the ACOPF (OPF
  Paper 4)},'' FERC, Tech. Rep., Mar. 2013.

\bibitem{molzahn_lesieutre_demarco-global_optimality_condition}
D.~K. Molzahn, B.~C. Lesieutre, and C.~L. DeMarco, ``{A Sufficient Condition
  for Global Optimality of Solutions to the Optimal Power Flow Problem},''
  \emph{IEEE Trans. Power Syst.}, vol.~29, no.~2, pp. 978--979, 2014.

\bibitem{ferc5}
A.~Castillo and R.~O'Neill, ``{Computational Performance of Solution Techniques
  Applied to the ACOPF (OPF Paper 5)},'' FERC, Tech. Rep., Jan. 2013.

\bibitem{molzahn_holzer_lesieutre_demarco-large_scale_sdp_opf}
D.~K. Molzahn, J.~T. Holzer, B.~C. Lesieutre, and C.~L. DeMarco,
  ``{Implementation of a Large-Scale Optimal Power Flow Solver Based on
  Semidefinite Programming},'' \emph{IEEE Trans. Power Syst.}, vol.~28, no.~4,
  pp. 3987--3998, 2013.

\bibitem{pscc2014}
D.~K. Molzahn and I.~A. Hiskens, ``{Moment-Based Relaxation of the Optimal
  Power Flow Problem},'' \emph{18th Power Syst. Comput. Conf. (PSCC)}, 18-22
  Aug. 2014.

\bibitem{madani2014}
R.~Madani, S.~Sojoudi, and J.~Lavaei, ``{Convex Relaxation for Optimal Power
  Flow Problem: Mesh Networks},'' \emph{IEEE Trans. Power Syst.}, vol.~30,
  no.~1, pp. 199--211, Jan. 2015.

\bibitem{molzahn_hiskens-sparse_moment_opf}
D.~Molzahn and I.~Hiskens, ``{Sparsity-Exploiting Moment-Based Relaxations of
  the Optimal Power Flow Problem},'' \emph{IEEE Trans. Power Syst.}, vol.~30,
  no.~6, pp. 3168--3180, Nov. 2015.

\bibitem{ibm_opf}
B.~Ghaddar, J.~Marecek, and M.~Mevissen, ``{Optimal Power Flow as a Polynomial
  Optimization Problem},'' \emph{IEEE Trans. Power Syst.}, vol.~31, no.~1, pp.
  539--546, Jan. 2016.

\bibitem{josz_molzahn-complex_hierarchy}
C.~Josz and D.~Molzahn, ``{Moment/Sums-of-Squares Hierarchy for Complex
  Polynomial Optimization},'' \emph{arXiv:1508.02068}, 2015.

\bibitem{cedric_msdp}
C.~Josz, J.~Maeght, P.~Panciatici, and J.~C. Gilbert, ``{Application of the
  Moment-SOS Approach to Global Optimization of the OPF Problem},'' \emph{IEEE
  Trans. Power Syst.}, vol.~30, no.~1, pp. 463--470, Jan. 2015.

\bibitem{low_tutorial}
S.~Low, ``{Convex Relaxation of Optimal Power Flow: Parts I \& II},''
  \emph{{IEEE Trans. Control Network Syst.}}, vol.~1, no.~1, pp. 15--27, Mar.
  2014.

\bibitem{coffrin2015qc}
C.~Coffrin, H.~Hijazi, and P.~{Van Hentenryck}, ``{The QC Relaxation: A
  Theoretical and Computational Study on Optimal Power Flow},'' \emph{IEEE
  Trans. Power Syst.}, vol.~31, no.~4, pp. 3008--3018, July 2016.

\bibitem{sun2015}
B.~Kocuk, S.~Dey, and A.~Sun, ``{Strong SOCP Relaxations of the Optimal Power
  Flow Problem},'' \emph{{\rm To appear in} Oper. Res.}, 2016.

\bibitem{ghaddar2015}
X.~Kuang, L.~F. Zuluaga, B.~Ghaddar, and J.~Naoum-Sawaya, ``{Approximating the
  ACOPF Problem with a Hierarchy of SOCP Problems},'' in \emph{IEEE PES General
  Meeting}, July 2015, pp. 1--5.

\bibitem{bienstock2014}
D.~Bienstock and G.~Munoz, ``{On Linear Relaxations of OPF Problems},''
  \emph{arXiv: 1411.1120}, Nov. 2014.

\bibitem{coffrin2016pscc}
C.~Coffrin, H.~Hijazi, and P.~{Van Hentenryck}, ``{Network Flow and Copper
  Plate Relaxations for AC Transmission Systems},'' \emph{19th Power Syst.
  Comput. Conf. (PSCC)}, 20-24 June 2016.

\bibitem{hiskens2001}
I.~Hiskens and R.~Davy, ``{Exploring the Power Flow Solution Space Boundary},''
  \emph{IEEE Trans. Power Syst.}, vol.~16, no.~3, pp. 389--395, 2001.

\bibitem{bernie_opfconvexity}
B.~Lesieutre and I.~Hiskens, ``{Convexity of the Set of Feasible Injections and
  Revenue Adequacy in FTR Markets},'' \emph{IEEE Trans. Power Syst.}, vol.~20,
  no.~4, pp. 1790--1798, Nov. 2005.

\bibitem{hill2008}
Y.~V. Makarov, Z.~Y. Dong, and D.~J. Hill, ``{On Convexity of Power Flow
  Feasibility Boundary},'' \emph{IEEE Trans. Power Syst.}, vol.~23, no.~2, pp.
  811--813, May 2008.

\bibitem{lavaei_geometry}
J.~Lavaei, D.~Tse, and B.~Zhang, ``{Geometry of Power Flows and Optimization in
  Distribution Networks},'' \emph{IEEE Trans. Power Syst.}, vol.~29, no.~2, pp.
  572--583, Mar. 2014.

\bibitem{chandra2015equilibria}
S.~Chandra, D.~Mehta, and A.~Chakrabortty, ``{Equilibria Analysis of Power
  Systems Using a Numerical Homotopy Method},'' in \emph{IEEE PES General
  Meeting}, July 26-30, 2015, pp. 1--5.

\bibitem{dj2015}
K.~Dvijotham and K.~Turitsyn, ``{Construction of Power Flow Feasibility
  Sets},'' \emph{arXiv:1506.07191}, 2015.

\bibitem{alves2002}
D.~Alves and G.~{da Costa}, ``{An Analytical Solutlon to the Optimal Power
  Flow},'' \emph{IEEE Power Eng. Rev.}, vol.~22, no.~3, pp. 49--51, Mar. 2002.

\bibitem{venikov1975}
V.~Venikov, V.~Stroev, V.~Idelchick, and V.~Tarasov, ``{Estimation of
  Electrical Power System Steady-State Stability in Load Flow Calculations},''
  \emph{IEEE Trans. Power App. Syst.}, vol.~94, no.~3, pp. 1034--1041, May
  1975.

\bibitem{voltage_collapse}
H.-D. Chiang, I.~Dobson, R.~Thomas, J.~Thorp, and L.~Fekih-Ahmed, ``{On Voltage
  Collapse in Electric Power Systems},'' \emph{IEEE Trans. Power Syst.},
  vol.~5, no.~2, pp. 601--611, May 1990.

\bibitem{cpf}
V.~Ajjarapu and C.~Christy, ``{The Continuation Power Flow: A Tool for Steady
  State Voltage Stability Analysis},'' \emph{IEEE Trans. Power Syst.}, vol.~7,
  no.~1, pp. 416--423, Feb. 1992.

\bibitem{alvarado1994}
F.~Alvarado, I.~Dobson, and Y.~Hu, ``{Computation of Closest Bifurcations in
  Power Systems},'' \emph{IEEE Trans. Power Syst.}, vol.~9, no.~2, pp.
  918--928, May 1994.

\bibitem{guo1990}
S.~Guo and F.~Salam, ``{Determining the Solutions of the Load Flow of Power
  Systems: Theoretical Results and Computer Implementation},'' in \emph{IEEE
  29th Ann. Conf. Decis. Contr. (CDC)}, Dec. 1990, pp. 1561--1566.

\bibitem{mehta2014a}
D.~Mehta, H.~Nguyen, and K.~Turitsyn, ``{Numerical Polynomial Homotopy
  Continuation Method to Locate All The Power Flow Solutions},''
  \emph{arXiv:1408.2732}, 2014.

\bibitem{SW:05}
A.~J. Sommese and C.~W. Wampler, \emph{{The Numerical Solution of Systems of
  Polynomials Arising in Engineering and Science}}.\hskip 1em plus 0.5em minus
  0.4em\relax World Scientific Publishing Company, 2005.

\bibitem{BHSW06}
D.~J. Bates, J.~D. Hauenstein, A.~J. Sommese, and C.~W. Wampler, ``{Bertini:
  Software for Numerical Algebraic Geometry},'' Available at
  http://www.nd.edu/$\sim$sommese/bertini.

\bibitem{paramotopy}
D.~Bates, D.~Brake, and M.~Niemerg, ``{Paramotopy: Parameter Homotopies in
  Parallel},'' Available at http://www.paramotopy.com.

\bibitem{illustrative_example}
D.~Molzahn and I.~Hiskens, ``{Convex Relaxations of Optimal Power Flow
  Problems: An Illustrative Example},'' \emph{IEEE Trans. Circuits Syst. I:
  Reg. Papers}, vol.~63, no.~5, pp. 650--660, May 2016.

\bibitem{acc2016}
D.~Molzahn, D.~Mehta, and M.~Niemerg, ``{Toward Topologically Based Upper
  Bounds on the Number of Power Flow Solutions},'' \emph{American Control Conf.
  (ACC)}, July 2016.

\bibitem{chen2015bounds}
T.~Chen and D.~Mehta, ``{On the Network Topology Dependent Solution Count of
  the Algebraic Load Flow Equations},'' \emph{arXiv:1512.04987}, 2015.

\bibitem{lasserre2001}
J.-B. Lasserre, ``{Global Optimization with Polynomials and the Problem of
  Moments},'' \emph{SIAM J. Optimiz.}, vol.~11, no.~3, pp. 796--817, 2001.

\bibitem{ibm_paper}
B.~Ghaddar, J.~Marecek, and M.~Mevissen, ``{Optimal Power Flow as a Polynomial
  Optimization Problem},'' \emph{IEEE Trans. Power Syst.}, vol.~31, no.~1, pp.
  539--546, Jan. 2016.

\bibitem{nie2014}
J.~Nie, ``{Optimality Conditions and Finite Convergence of Lasserre's
  Hierarchy},'' \emph{Math. Program.}, vol. 146, no. 1-2, pp. 97--121, 2014.

\bibitem{coffrin_tightening}
C.~Coffrin, H.~Hijazi, and P.~{Van Hentenryck}, ``{Strengthening Convex
  Relaxations with Bound Tightening for Power Network Optimization},'' in
  \emph{Principles and Practice of Constraint Programming}, ser. Lecture Notes
  in Computer Science, G.~Pesant, Ed.\hskip 1em plus 0.5em minus 0.4em\relax
  Springer International Publishing, 2015, vol. 9255, pp. 39--57.

\bibitem{chen2015}
C.~Chen, A.~Atamt\"urk, and S.~Oren, ``{Bound Tightening for the Alternating
  Current Optimal Power Flow Problem},'' \emph{{\rm to appear in} IEEE Trans.
  Power Syst.}

\bibitem{yalmip}
J.~Lofberg, ``{YALMIP: A Toolbox for Modeling and Optimization in MATLAB},'' in
  \emph{{IEEE Int. Symp. Compu. Aided Control Syst. Des.}}, 2004, pp. 284--289.

\bibitem{bertinilab}
D.~J. Bates, A.~J. Newell, and M.~Niemerg, ``{{BertiniLab: A MATLAB Interface
  for Solving Systems of Polynomial Equations}},'' \emph{Numer. Algorithms},
  vol.~71, no.~1, pp. 229--244, 2015.

\bibitem{hicss2014}
D.~K. Molzahn, B.~C. Lesieutre, and C.~L. DeMarco, ``{Investigation of Non-Zero
  Duality Gap Solutions to a Semidefinite Relaxation of the Power Flow
  Equations},'' in \emph{47th Hawaii Int. Conf. Syst. Sci. (HICSS)}, 6-9 Jan.
  2014.

\end{thebibliography}
%

%

\vspace*{-3em}
\begin{IEEEbiography}[{\includegraphics[width=1in,height=1.25in,clip,keepaspectratio]{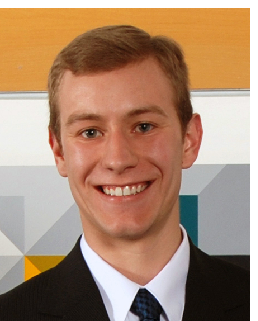}}]{Daniel K. Molzahn}
(S'09-M'13) is a Computational Engineer at Argonne National Laboratory. He was a Dow Postdoctoral Fellow in Sustainability at the University of Michigan, Ann Arbor and received the B.S., M.S., and Ph.D. degrees in electrical engineering and the Masters of Public Affairs degree from the University of Wisconsin--Madison, where he was a National Science Foundation Graduate Research Fellow. His research focuses on optimization and control of electric power systems.
\end{IEEEbiography}









\end{document}